\documentclass[conference]{IEEEtran}

\usepackage{amsmath}
\usepackage{amssymb}
\usepackage{graphicx}
\usepackage{balance}

\def\<{\leqslant}           
\def\>{\geqslant}           

\def\[{[\![}
\def\]{]\!]}

\newcommand{\bx}{\bar{x}}
\newcommand{\bt}{\bar{t}}
\newcommand{\bD}{\bar{D}}
\newcommand{\bA}{\bar{A}}
\newcommand{\bG}{\bar{G}}
\newcommand{\bC}{\bar{C}}

\newcommand{\bPi}{\bar{\Pi}}
\newcommand{\hx}{\hat{x}}
\newcommand{\tx}{\tilde{x}}
\newcommand{\ha}{\hat{a}}
\newcommand{\hc}{\hat{c}}
\newcommand{\hg}{\hat{g}}
\newcommand{\hal}{\hat{\alpha}}
\newcommand{\hka}{\hat{\kappa}}
\newcommand{\hxi}{\hat{\xi}}
\newcommand{\tV}{\tilde{V}}

\begin{document}
\title{Robust Smoothing for Discrete-Time Uncertain Nonlinear Systems}

\author{\IEEEauthorblockN{Abhijit G. Kallapur}
\IEEEauthorblockA{School of Engineering and Information Technology\\
University of New South Wales\\
Australian Defence Force Academy\\
Canberra, ACT 2600, Australia\\
Email: a.kallapur@adfa.edu.au}

\and
\IEEEauthorblockN{Ian R. Petersen}
\IEEEauthorblockN{School of Engineering and Information Technology\\
University of New South Wales\\
Australian Defence Force Academy\\
Canberra, ACT 2600, Australia\\
Email: i.petersen@adfa.edu.au}}

\maketitle


\begin{abstract}
This paper derives recursion equations for a robust smoothing problem for a class of nonlinear systems with uncertainties in modeling and exogenous noise sources. The systems considered operate in discrete-time and the uncertainties are modeled in terms of a sum quadratic constraint. The robust smoothing problem is solved in terms of a forward-time and a reverse-time filter. Both these filters are formulated in terms of set-valued state estimators and are recast into subsidiary optimal control problems. These optimal control problems are described in terms of discrete-time Hamilton-Jacobi-Bellman equations, whose approximate solutions lead to recursive Riccati difference equations, filter state equations, and level shift scalar equations for the forward-time and the reverse-time filters.  
\end{abstract}

\IEEEpeerreviewmaketitle

\section{Introduction}
Filtering and smoothing are extensively used for estimation of data from noisy signals. Although filtering by itself is sufficient in most cases, smoothing becomes necessary where outliers need to be estimated from noisy data. These topics have been of interest in the estimation and control community since the 1960s; see e.g., \cite{kalman1960,meditch1969,anderson1979,bar-shalom1993}. In spite of the advancements in the theory for estimation and smoothing, the success of most of the algorithms are based on precise knowledge of the system dynamics. In addition, such methods are known to diverge in the presence of unmodeled system dynamics and parametric uncertainties \cite{theodor1994}. In order to take into account uncertainties in modeling, robust versions of these estimators were developed in the stochastic as well as deterministic settings; see e.g., \cite{savkin1995,james1998,petersen1999,petersen2000,wang2002,yoon2004,dong2006,souto2009,ieeetac2009,acc2009,ifac2011,acc2012}. Also, the issue of robust smoothing for systems with model uncertainties have been addressed in the stochastic and $H_{\infty}$ frameworks in \cite{martin1983,blanco2001,malcom2005,felsberg2006,hongguo2007,garcia2010}.

As seen from research in the past, most of the robust filtering, prediction, and smoothing techniques are not applicable to systems where the statistical nature of noise is unclear and hence an appropriate statistical model for the system is not available. Such a situation may arise in systems where the uncertainties in modeling dominate over uncertainties due to exogenous noise sources. In such cases, a deterministic approach to filtering and smoothing might be more appropriate. Although there are various approaches to deterministic filtering, \cite{bertsekas1971} provides a set-membership approach to Kalman filtering. This set-membership approach to deterministic filtering involves computing the set of all possible states compatible with given output measurements for an uncertain system with norm bound uncertainties. This method was extended by Savkin and Petersen in \cite{savkin1995,petersen1999} to linear uncertain systems where the uncertainties were bound by integral and sum quadratic constraints. Furthermore, the method in \cite{savkin1995,petersen1999} was extended to include nonlinear uncertain systems in \cite{james1998,ieeetac2009,acc2009,ifac2011,acc2012}. However, all these methods deal with the problem of deterministic filtering for uncertain systems and do not deal with issues related to robust smoothing. The set-membership state estimation approach in \cite{savkin1995,petersen1999} was used to derive a robust smoothing algorithm for uncertain linear systems in \cite{moheimani1998}. This method was extended to uncertain nonlinear systems in the continuous-time case in \cite{acc2013}.

In this paper, we present a robust smoothing algorithm for uncertain nonlinear systems in a discrete-time setting, with uncertainties bound by a sum quadratic constraint (SQC). Such a smoothing algorithm is applicable for systems in fields such as aerospace engineering and signal processing where discrete-time sensors are used. In this paper, we derive a fixed-point smoothing algorithm that solves a bridge filtering problem. This bridging approach to smoothing comprises of solving two underlying robust filtering problems. The corresponding solution leads to a forward-time filter over the discrete-time interval $[0,k]$ and a reverse-time filter over the discrete-time interval $[k,t]$. Here, $0 \le k \le t$.

Both the forward-time and reverse-time filtering problems are formulated as set-valued state estimation problems which are then recast into optimal control problems. The optimal control problems are described in terms of Hamilton-Jacobi-Bellman (HJB) equations evaluated over different time intervals. Various nonlinear systems in the HJB equations are linearized and approximate quadratic solutions are computed for the dynamic programming equations. This leads to a Riccati equation, a filter state equation, and a level shift scalar equation for both the forward-time and the reverse-time filters. Although it is straightforward to formulate the set-valued state estimation problem for the reverse-time filter, this is not trivial in the case of the forward-time filter. The set-valued state estimation problem for the forward-time filter will be solved by obtaining an approximate solution to the optimal control problem in terms of a forward-time dynamic programming equation for a reverse-time discrete-time uncertain nonlinear system with uncertainties described in terms of an SQC. For a detailed explanation, see \cite{thesis2009}.

The rest of the paper is organized as follows: Section \ref{sec:problem-formulation} introduces the formulation of the discrete-time uncertain nonlinear system starting from an equivalent continuous-time system. The set-valued state estimation approach to deterministic filtering is presented in Section \ref{sec:svse}. The formulation of the robust smoothing problem in terms of a forward-time and a reverse-time filter are presented in Section \ref{sec:robust-smoothing-problem}. The recursion equations for the discrete-time robust nonlinear smoother in terms of the Riccati equation, the filter state equation, and the level shift scalar equation are derived in Section \ref{sec:robust-smoother}. The paper concludes with a summary and future research directions in Section \ref{sec:conclusion}.
\section{Problem Formulation}
\label{sec:problem-formulation}
In this section, we introduce a general model for a continuous-time nonlinear uncertain system where the uncertainties are modeled in terms of an integral quadratic constraint (IQC). Furthermore, we provide a discussion on obtaining an equivalent discrete-time uncertain system where the uncertainties are modeled in terms of an SQC.

\begin{figure}[h]
\centering{\includegraphics[width=8.5cm]{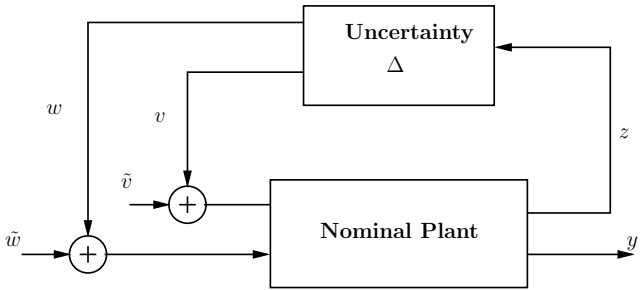}}
\caption{Uncertain system.}
\label{fig:uncertain-sys}
\end{figure}

Consider the following forward-time continuous-time uncertain nonlinear system
\begin{align}
	&\dot{x}(s) = a_c(s, x(s), u(s)) + D_c w(s), \nonumber \\
	&z(s) = g_c(s, x(s), u(s)), \nonumber \\
\label{eq:ct-sys}
	&y(s) = c_c(s,x(s),u(s)) + v(s).
\end{align}
Here, $x(s) \in \mathbb{R}^n$ is the state vector, $z(s) \in \mathbb{R}^q$ is the uncertainty output, $y(s) \in \mathbb{R}^l$ is the measured output, $w(s) \in \mathbb{R}^p$ and $v(s) \in \mathbb{R}^l$ are the process and measurement uncertainty inputs. Also, $a_c(\cdot), g_c(\cdot), \text{and} c_c(\cdot)$ are given sufficiently smooth nonlinear functions of the state variables and $D_c$ is a matrix function of time. Such a system is depicted in Fig. \ref{fig:uncertain-sys}.

The uncertainty inputs $w(s) \;\text{and} \; v(s)$ are assumed to be square integrable, satisfying an IQC at any time $t \ge 0$ as in \cite{james1998,acc2012},
\begin{align}
	S_c&(t,x(t),w([0,t]) 
	:= 
	\| x(0) - \bx_0 \|_N^2 \nonumber \\
	&+ \int_0^t
	\left(
		\|w(s)\|_{Q(s)}^2 + \|v(s)\|_{R(s)}^2 - |z(s)|^2
	\right)
	ds
	< d.
\label{eq:iqc}	
\end{align}
Here, $x_0$ is the initial state value, $\bx_0$ is the nominal initial state, and $w_{[0,t]} := w(s)|_{0 \le s \le t}$ represents the restriction of the noise term $w(\cdot)$ over the time interval $[0,t]$. All vectors in \eqref{eq:iqc} and in the rest of the paper are organized as columns, unless specified otherwise. Also, $w(s)$ and $v(s)$ are admissible uncertainties described by
\begin{equation}
	\begin{bmatrix}
		w(s) \\
		v(s)
	\end{bmatrix}
	=
	\eta(s,x(\cdot)),
\label{eq:uncertainty}
\end{equation}
where $\eta(\cdot)$ is a nonlinear time varying uncertain function and $x(\cdot)$ represents the state vector $x(s)$ for all instants of time $s \in [0,t]$. The fact that the uncertainty inputs $w(s)$ and $v(s)$ depend on $x(\cdot)$ in \eqref{eq:uncertainty} indicates that they are allowed to depend dynamically on the system state. Also, $N = N^T \succ 0$ is a given matrix and $Q(s) = Q(s)^T \succ 0, R(s) = R(s)^T \succ 0$ are given matrix functions of time.

In order to derive recursion equations for the discrete-time robust smoother, we need to convert the continuous-time nonlinear system in \eqref{eq:ct-sys} and the corresponding IQC \eqref{eq:iqc} into corresponding discrete forms. As mentioned in the introduction, the recursion equations for the smoother will be obtained by solving for a forward-time filter and a reverse-time filter. The derivation of the forward-time filter in the set-valued state estimation framework, is most straightforward if the continuous-time uncertain system in \eqref{eq:ct-sys} is discretized in reverse-time rather than in forward-time. The application of higher-order Euler or Runge-Kutta methods to \eqref{eq:ct-sys} gives us the following reverse-time discrete-time system
\begin{align}
	&x_t = a_t(x_{t+1}) - D_t w_t, \nonumber \\
	&z_t = g_t(x_{t+1}), \nonumber \\
\label{eq:rtdt-sys}
	&y_t = c_t(x_{t+1}), + v_{t+1},
\end{align}
where $a_t, g_t, c_t$ are discrete-time nonlinear functions and $D_t$ is a given matrix function of time, all expressed at time step $t$.

Indeed, the solution to the reverse-time filter in the set-valued state estimation framework is obtained by discretizing the continuous-time uncertain system \eqref{eq:ct-sys} forward in time. Applying higher-order Euler or Runge-Kutta methods to \eqref{eq:ct-sys}, we arrive at the following forward-time discrete-time uncertain nonlinear system
\begin{align}
	&x_{t+1} = \alpha_t(x_t) + \bar{D}_t w_t, \nonumber \\
	&z_{t+1} = \kappa_t(x_{t+1}), \nonumber \\
\label{eq:ftdt-sys}
	&y_{t+1} = \xi_t(x_{t+1}), + v_{t+1},
\end{align}
where $\alpha_t, \kappa_t, \xi_t$ are given nonlinear functions of time and $\bar{D}_t$ is a given time-varying matrix at time step $t$. The nonlinear functions in \eqref{eq:rtdt-sys} and \eqref{eq:ftdt-sys} include the effect of known control signal $u_t$ at time step $t$. However, their explicit dependence on the control signal $u_t$ will be omitted as the paper deals with a smoothing problem.

In the case of the forward-time filter and the reverse-time filter, the uncertainties associated with the systems \eqref{eq:rtdt-sys} and \eqref{eq:ftdt-sys} are modeled in terms of an SQC. This SQC is derived by discretizing the IQC in \eqref{eq:iqc} using Reimann sum approximation to give
\begin{align}
	S_t&(x_t, w_{[0,t]}) \nonumber \\
	&:= 
	\| x_0 - \bx_0 \|_N^2 
	+ \sum_{s=0}^{t-1}
	\left(
		\|w_s\|_{Q_s}^2 + \|v_s\|_{R_s}^2 - |z_s|^2
	\right)
	< d.
\label{eq:sqc}	
\end{align}
Here, $N = N^T \succ 0$ is a given matrix and $Q_s = Q_s^T \succ 0, R_s = R_s^T \succ 0$ are given matrix functions of time. 

Since the recursion equations for the robust smoother are solved in terms of a forward-time filter and a reverse-time filter in the set-valued state estimation framework, the following section introduces the concept of set-valued state estimation.

\section{Set-valued State Estimation}
\label{sec:svse}
In the continuous-time case, consider an output $y_0(\cdot)$ for any $s \in [0,t]$ for the system in \eqref{eq:ct-sys}. The set-valued state estimation problem comprises of finding the set 
\begin{equation}
\label{eq:ct-svse}
	Z_s\left[ \bx_0, y_0(\cdot)|_0^s, d \right]
\end{equation}
of all possible states $x(s)$ at time $s$ for the system \eqref{eq:ct-sys} with initial conditions and uncertainties defined by the IQC \eqref{eq:iqc}, consistent with the output $y_0(\cdot)$. Given an output sequence $y_0(\cdot)$, it follows from the definition of $Z_s$ that
\begin{equation}
\label{eq:svse}
	x_s
	\in
	Z_s 
	\left[
		\bx_0, y_0(\cdot)|_0^s, d
	\right]
\end{equation} 
if and only if there exists an uncertain input sequence $w(\cdot)$ such that, $J[x_s, w(\cdot)] \le d$, where the cost functional $J[x_s, w(\cdot)]$ is defined in terms of the IQC \eqref{eq:iqc} as
\begin{equation}
\label{eq:ct-oc}
	J[x_s, w(\cdot)]
	\triangleq
	\inf_{w(\cdot)} S_c(t,x(t),w(\cdot))]
\end{equation}
with $y_0(t) = c_c(x(t)) + v(t), z(t) = g_c(x(t))$. Now we present the formulation of the discrete-time robust smoothing problem in the set-valued state estimation framework.

\section{Robust Smoothing and the Optimal Control Problem}
\label{sec:robust-smoothing-problem}
In this section, we consider a discrete-time robust smoothing problem formulated in the set-valued state estimation framework. Here, the robust smoothing problem is considered at a fixed point in time $k = (t-\bt)$ and comprises of computing all the states $x_k \in X_k[\bx_0, y_{[0,k]}, d]$ at time $k$. Indeed, the corresponding optimal control problem in this case, will depend on the time interval being considered. 

From the definition for the set-valued state estimator in \eqref{eq:ct-svse}-\eqref{eq:ct-oc}, at time $k$ we have in the discrete-time case
\begin{align}
	&X_{k}[\bx_0, y_{[0,t]},d]
	= 
	\Big\{
		x_{k} \in \mathbb{R}^n
		:
		\inf_{w_{[0,t]}}
		\Big(
		S_1(x_{k}, w_{[0,k]}) \nonumber \\
		& \qquad \qquad \qquad \qquad \qquad \qquad \qquad+
		S_2(x_{k}, w_{[k,t]})
		\Big)
		\le
		d
	\Big\},
\label{eq:Xtmbt}
\end{align}
where the infimum is taken over all $x_{(\cdot)}$ and $w_{[0,t]}$ with the boundary condition $x_{k} = \bx_{k}$. The value function $S_1(\cdot)$ is defined over the time interval $[0,k]$ as
\begin{align}
	S_1&(x_{k}, w_{[0,k]}) \nonumber \\
	&:= 
	\| x_0 - \bx_0 \|_N^2 
	+ \sum_{s=0}^{k}
	\left(
		\|w_s\|_{Q_s}^2 + \|v_s\|_{R_s}^2 - |z_s|^2
	\right).
\label{eq:S1}
\end{align}
The value function $S_2(\cdot)$ is defined over the time interval $[k,t]$ as
\begin{equation}
	S_2(x_t, w_{[k,t]}) 
	\!:=\! 
	\sum_{s=k}^t
	\left(
		\|w_s\|_{Q_s}^2 + \|v_s\|_{R_s}^2 - |z_s|^2
	\right).
\label{eq:S2}
\end{equation}
Note the absence of the contribution of the error in initial state in \eqref{eq:S2}. This is because the initial state is known a priori over the time interval $[k, t]$. Indeed, for each $x_{k}$, $S_1(x_{k}, w_{[0,t]})$ depends on the sequence $w_{[0,t]}$ over the time interval $[0,k]$, whereas, $S_2(x_{k}, w_{[k,t]})$ depends on the sequence $w_{[0,t]}$ over the time interval $[k, t]$. Hence, the left hand side in the optimization problem in \eqref{eq:Xtmbt} can be split into the sum of two optimization problems based on the time interval being considered. The first optimization problem depends on the infimum computed over $w_{(\cdot)} \in \mathcal{L}^2[0,k]$, whereas the second one depends on the infimum computed over $w_{(\cdot)} \in \mathcal{L}^2[k,t]$. Hence, the definition for $X_{k}[\bx_0, y_{[0,t]},d]$ in \eqref{eq:Xtmbt} can be restated as
\begin{align}
	&X_{k}[\bx_0, y_{[0,t]},d] 
	= 
	\Big\{
		x_{k} \in \mathbb{R}^n
		:
		\inf_{w_{(\cdot)}\in \mathcal{L}^2[0,k]}
		\Big( S_1(x_{k}, w_{(\cdot)}) \Big) \nonumber \\
		& \qquad \qquad \qquad \qquad +
		\inf_{w_{(\cdot)}\in \mathcal{L}^2[k,t]}
		\Big( S_2(x_{k}, w_{(\cdot)}) \Big)
		\le
		d.
	\Big\}
\label{eq:Xtmbt-split}
\end{align}

The optimization problem 
\begin{equation}
	\inf_{w_{(\cdot)} \in \mathcal{L}^2[0,k]}
		S_1(x_{k}, w_{(\cdot)})
	+
	\inf_{w_{(\cdot)}\in \mathcal{L}^2[k,t]}
		S_2(x_{k}, w_{(\cdot)})
\label{eq:opt-prob}
\end{equation}
for the systems \eqref{eq:rtdt-sys} and \eqref{eq:ftdt-sys} satisfying the quadratic constraints defined in \eqref{eq:S1} and \eqref{eq:S2} respectively, defines a nonlinear optimal control problem with a sign indefinite quadratic cost function. The robust smoothing problem will be solved by obtaining an approximate solution to this optimal control problem.

\section{Robust Nonlinear Smoother}
\label{sec:robust-smoother}
In order to obtain recursion equations for the robust nonlinear smoothing problem, the two optimization problems in \eqref{eq:opt-prob} will be individually solved using methods of dynamic programming. This is accomplished by finding approximate solutions to the underlying HJB equations. Note that the solution to the robust nonlinear smoothing problem consists of two filters: one operating in forward-time over the interval $[0,k]$ and the other operating in reverse-time over the interval $[k,t]$. The solution to the forward-time filter will comprise of approximately solving the HJB equation forward in time for the reverse-time discrete-time uncertain nonlinear system in \eqref{eq:rtdt-sys} and the uncertainties described by the SQC \eqref{eq:S1}. The solution to the reverse-time filter will involve solving for the HJB equation reverse in time with the forward-time discrete-time uncertain nonlinear system \eqref{eq:ftdt-sys} and the uncertainties described by the SQC \eqref{eq:S2}. The formulations and solutions to the optimal control problems for the forward-time and reverse-time filters will be discussed in Sections \ref{ssec:ft-filter} and \ref{ssec:rt-filter} respectively.

\subsection{Forward-Time Filter}
\label{ssec:ft-filter}
In order to derive the Riccati equation, the filter equation, and the level shift scalar equation for the forward-time filter, we consider the following optimization problem involving $S_1(\cdot)$ from \eqref{eq:opt-prob} defined over the time interval $[0,k]$,
\begin{equation}
	H_k(x) 
	:= 
	\inf_{w_{(\cdot)} \in \mathcal{L}^2[0,k]}
		S_1(x_{k}, w_{(\cdot)}).
\label{eq:ftf-opt}
\end{equation}
Here, $H_k(x)$ defines a reverse-time optimal control problem with $S_1(\cdot)$ defined in \eqref{eq:S1}. In the continuous-time case, this optimization problem can be solved using dynamic programming methods in terms of the forward-time system dynamics. However, this is not straightforward in the discrete-time case, since it would involve inverting the nonlinear function describing the discrete-time nonlinear system, which may not be possible. In order to solve for the set-valued state estimation problem in this case, we consider an approximate solution for the optimal control problem \eqref{eq:ftf-opt} in terms of a forward-time dynamic programming equation for the reverse-time discrete-time uncertain nonlinear system \eqref{eq:rtdt-sys} with an SQC uncertainty description \eqref{eq:S1}. This technique was first proposed in \cite{ieeetac2009,acc2009}.

The corresponding HJB equation for the optimization problem in \eqref{eq:ftf-opt} for the reverse-time discrete-time system \eqref{eq:rtdt-sys} and the SQC \eqref{eq:S1} is given by
\begin{align}
	V_{k+1}(x, y_{0:k}) 
	= 
	\inf_{w_{(\cdot)} \in \mathcal{L}^2[0,k]}
	&\big(
		V_k 
		\left(
			a_k(x) - D_k w_k
		\right) \nonumber \\
\label{eq:ftf-hjb}
	&+
	\|w_k\|_{Q_s}^2 
	+ 
	\|v_k\|_{R_s}^2 
	- 
	|z_k|^2
	\big)
\end{align}
with initial condition
\begin{equation}
	V_0(x) = \| x - \bx_0 \|^2_N.
\end{equation}

As a first step towards obtaining an approximate solution to the HJB \eqref{eq:ftf-hjb}, we approximate the value function $V_k(\cdot)$ as
\begin{equation}
\label{eq:ftf-quad-approx}
	V_k(x) 
	\approx
	\|
		x - \hx_k 
	\|^2_{\Pi_k}
	+
	\phi_k.
\end{equation}
Applying the approximate solution \eqref{eq:ftf-quad-approx} to the HJB \eqref{eq:ftf-hjb} we get
\begin{align}
	\|
		x - \hx_k 
	\|^2_{\Pi_{k+1}}
	+
	&\phi_{k+1}
	=
	\inf_{w_{(\cdot)} \in \mathcal{L}^2[0,k]}
	\big(
		\|
			a_k(x) - D_k w_k
		\|^2_{\Pi_k} \nonumber \\
\label{eq:ftf-1}
	&+
		\|w_k\|_{Q_k}^2 
		+ 
		\|v_k\|_{R_k}^2 
		- 
		|z_k|^2
		+
		\phi_k
	\big).
\end{align}
Solving for the infimum problem on the right hand side of \eqref{eq:ftf-1} using the method of completing the squares, we obtain the following optimum value for $w_k$, denoted as $w^*$
\begin{equation}
\label{eq:ftf-w*}
	w^*
	=
	(
		D_k^T \Pi_k D_k + Q_k
	)^{-1}
	D_k^T \Pi_k
	(
		a_k(x) - \hx_k
	).
\end{equation}
Substituting $w_k = w^*$ from \eqref{eq:ftf-w*} into \eqref{eq:ftf-1} and rearranging terms we get
\begin{align}
	\|
		x - \hx_k 
	\|^2_{\Pi_{k+1}}
	+
	&\phi_{k+1}
	=
	\|
		a_k(x)
	\|^2_{\Sigma}
	-
	2 \hx_k^T \Sigma a_k(x) \nonumber \\
\label{eq:ftf-2}
	&+
	\|
		\hx_k
	\|^2_{\Sigma}
	+
	\|
		y_k - c_k(x)
	\|^2_{R_k}
	-
	|g(x)|^2
	+
	\phi_k,
\end{align}
where
\begin{equation}
\label{eq:Sigma}
	\Sigma 
	= 
	\Pi_k
	-
	\Pi_k D_k
	(
		D_k^T \Pi_k D_k + Q_k
	)^{\#}
	D_k^T \Pi_k.
\end{equation}
Here, $(\cdot)^{\#}$ denotes the Moore-Penrose pseudo-inverse. If the matrix $(D_k^T \Pi_k D_k + Q_k)$ is positive-definite, the pseudo-inverse in \eqref{eq:Sigma} can be replaced by the normal matrix inverse. This is feasible and holds in most cases where a suitable solution exists for $\Pi_k$; see e.g., \cite{petersen1999}.

In order to obtain an approximate solution to the right hand side of \eqref{eq:ftf-2}, we consider the following first order linearization for various nonlinear functions about the point $\hx_k$.
\begin{align}
	&a_k(x) \approx \ha_k + A_k (x - \hx_k), \nonumber \\
	&g_k(x) \approx \hg_k + G_k (x - \hx_k), \nonumber \\
\label{eq:ftf-lin}
	&c_k(x) \approx \hc_k + C_k (x - \hx_k),
\end{align}
where $\ha_k = a_k(\hx_k), \hg_k = g_k(\hx_k),$ and $\hc_k = c_k(\hx_k)$. Also, $A_k, G_k, $ and $C_k$ are the Jacobian matrices of the maps $a_k, g_k, $ and $c_k$ respectively.

Substituting for various linear approximations from \eqref{eq:ftf-lin} into \eqref{eq:ftf-2} we obtain
\begin{align}
	&\|
		x - \hx_k 
	\|^2_{\Pi_{k+1}}
	+
	\phi_{k+1}
	=
	\|
		x
	\|^2_{A_k^T \Sigma A_k + C_k^T \Sigma C_k - G_k^T G_k} \nonumber \\
	&\;\;-
	2 x^T
	\Big(
		\left(
			A_k^T \Sigma A_k
			+
			A_k^T \Sigma
			+
			C_k^T \Sigma C_k 
			- 
			G_k^T G_k
		\right) \hx_k \nonumber \\
		&\;\;\;\;\;\;\;\;\;\;-
		A_k^T \Sigma \ha_k
		+
		C_k^T R (y - \hc_k)
		+
		G_k^T G_k
	\Big) \nonumber \\
	&\;\;+
	\|
		\hx_k
	\|^2_{A_k^T \Sigma A_k + \Sigma - 2 \Sigma A_k + C_k^T R_k C_k - G_k^T G_k} \nonumber \\
	&\;\;\;-
	2 \hx_k^T
	\Big(
		A_k^T \Sigma \ha_k
		+
		\Sigma \ha_k
		-
		C_k^T R_k (y_k - \hc_k)
		-
		G_k^T \hg_k
	\Big) \nonumber \\
\label{eq:ftf-3}
	&\;\;\;+
	\ha_k^T \Sigma \ha_k
	+
	\| y_k \|^2_{R_k}
	-
	2 y_k^T R_k \hc_k
	+
	\| \hc_k \|^2_{R_k}
	-
	\hg_k^T \hg_k
	+
	\phi_k.
\end{align}
Comparing terms on the left hand side and right hand side of \eqref{eq:ftf-3}, we obtain the following recursion equations for the forward-time filter corresponding to the time interval $[0,k]$:
\begin{align}
	&\textit{Riccati equation over the time interval}\; [0,k]: \nonumber \\
	&\Pi_{k+1}
	=
	A_k^T \Sigma A_k 
	+ 
	C_k^T \Sigma C_k 
	- G_k^T G_k, \nonumber \\
\label{eq:ft-riccati}
	&\Pi_0 = N. \\
	&\textit{Filter state equation over the time interval}\; [0,k]: \nonumber \\
	&\hx_{k+1}
	=
	\hx_k
	+
	\Pi_k^{-1} \Xi, \nonumber \\
\label{eq:ft-filter}
	&\hx_0 = \bx_0, \\
	&\text{where} \nonumber \\
	&\Xi
	=
	\big(
		A_k^T \Sigma (\hx_k - \ha_k)
		+
		C_k^T R (y - \hc_k)
		+
		G_k^T G_k	
	\big). \nonumber \\
	&\textit{Level shift scalar equation over the time interval} \; [0, k]: \nonumber \\
	&\phi_{k+1}
	=
		\phi_k
		+
		\|
		\hx_k
	\|^2_{A_k^T \Sigma A_k + \Sigma - 2 \Sigma A_k + C_k^T R_k C_k - G_k^T G_k} \nonumber \\
	&\;\;\;-
	2 \hx_k^T
	\Big(
		A_k^T \Sigma \ha_k
		+
		\Sigma \ha_k
		-
		C_k^T R_k (y_k - \hc_k)
		-
		G_k^T \hg_k
	\Big) \nonumber \\
	&\;\;\;+
	\ha_k^T \Sigma \ha_k
	+
	\| y_k \|^2_{R_k}
	-
	2 y_k^T R_k \hc_k
	+
	\| \hc_k \|^2_{R_k}
	-
	\hg_k^T \hg_k, \nonumber \\
\label{eq:ft-scalar}
	&\phi_0 = 0.
\end{align}

\subsection{Reverse-Time Filter}
\label{ssec:rt-filter}
Here, we derive the Riccati equation, the filter state equation, and the level shift scalar equation for the reverse-time filter over the time interval $[k,t]$. This comprises of solving for the infimum involving the function $S_2(\cdot)$ in \eqref{eq:opt-prob}. This is given by
\begin{equation}
	B_k(x) 
	:= 
	\inf_{w_{(\cdot)} \in \mathcal{L}^2[k,t]}
		S_2(x_{k}, w_{(\cdot)}),
\label{eq:rtf-opt}
\end{equation}
where $S_2(\cdot)$ is defined in \eqref{eq:S2}. The corresponding reverse-time discrete-time HJB equation for the system \eqref{eq:ftdt-sys} and the SQC \eqref{eq:S2} over the time interval $[k,t]$ is given by
\begin{align}
	\tV_{k}(x, y_{0:k}) 
	= 
	\inf_{w_{(\cdot)} \in \mathcal{L}^2[k,t]}
	&\big(
		\tV_{k+1} 
		\left(
			\alpha_k(x) + \bD_k w_k
		\right) \nonumber \\
\label{eq:rtf-hjb}
	&+
	\|w_k\|_{Q_s}^2 
	+ 
	\|v_k\|_{R_s}^2 
	- 
	|z_k|^2
	\big)
\end{align}
with
\begin{equation}
	\tilde{V}_t(x) = 0.
\end{equation}

As mentioned in \ref{ssec:ft-filter} a first step towards obtaining an approximate solution to the HJB \eqref{eq:rtf-hjb}, we approximate the value function $\tilde{V}_k(\cdot)$ as
\begin{equation}
\label{eq:rtf-quad-approx}
	\tilde{V}_k(x) 
	\approx
	\|
		x - \tx_k 
	\|^2_{\bPi_k}
	+
	\psi_k.
\end{equation}
Applying the approximate solution \eqref{eq:rtf-quad-approx} to the HJB \eqref{eq:rtf-hjb} we get
\begin{align}
	\|
		x - \tx 
	\|^2_{\bPi_{k}}
	+
	\psi_{k}
	=
	&\inf_{w_{(\cdot)} \in \mathcal{L}^2[k,t]}
	\big(
		\|
			\alpha_k(x) + \bD_k w_k
		\|^2_{\bPi_{k+1}} 
		+ \psi_{k+1} \nonumber \\
\label{eq:rtf-1}
	&+
		\|w_k\|_{Q_k}^2 
		+ 
		\|v_k\|_{R_k}^2 
		- 
		|z_k|^2
	\big).
\end{align}
Solving for the infimum problem on the right hand side of \eqref{eq:rtf-1} using the method of completing the squares, we obtain the following optimum value for $w_k$, denoted as $w^*$
\begin{equation}
\label{eq:rtf-w*}
	w^*
	=
	(
		\bD_k^T \bPi_{k+1} \bD_k + Q_k
	)^{-1}
	\bD_k^T \bPi_{k+1}
	(
		\tx_k - \alpha_k(x)
	).
\end{equation}
Substituting $w_k = w^*$ from \eqref{eq:rtf-w*} into \eqref{eq:rtf-1} and rearranging terms we get
\begin{align}
	\|
		x - \tx_k 
	\|^2_{\bPi_{k}}
	+
	&\psi_{k}
	=
	\|
		\alpha_k(x)
	\|^2_{\Omega}
	-
	2 \tx_k^T \Omega \alpha_k(x) \nonumber \\
\label{eq:rtf-2}
	&+
	\|
		\tx_k
	\|^2_{\Omega}
	+
	\|
		y_k - \xi_k(x)
	\|^2_{R_k}
	-
	|\kappa(x)|^2
	+
	\psi_{k+1},
\end{align}
where
\begin{equation}
\label{eq:Omega}
	\Omega
	= 
	\bPi_{k+1}
	-
	\bPi_{k+1} \bD_k
	(
		\bD_k^T \bPi_{k+1} \bD_k + Q_k
	)^{\#}
	\bD_k^T \bPi_{k+1}.
\end{equation}
Here, $(\cdot)^{\#}$ denotes the Moore-Penrose pseudo-inverse. If the matrix $(\bD_k^T \bPi_k \bD_k + Q_k)$ is positive-definite, the pseudo-inverse in \eqref{eq:Omega} can be replaced by the normal matrix inverse. This is feasible and holds in most cases where a suitable solution exists for $\bPi_k$; see e.g., \cite{petersen1999}.

In order to obtain an approximate solution to the right hand side of \eqref{eq:rtf-2}, we consider first order linearization of various nonlinear functions in \eqref{eq:rtf-2} about the point $\tx_k$ as
\begin{align}
	&\alpha_k(x) \approx \hal_k + \bA_k (x - \tx_k), \nonumber \\
	&\kappa_k(x) \approx \hka_k + \bG_k (x - \tx_k), \nonumber \\
\label{eq:rtf-lin}
	&\xi_k(x) \approx \hxi_k + \bC_k (x - \tx_k),
\end{align}
where $\hat{\alpha}_k = \alpha_k(\tx_k), \hat{\kappa}_k = \kappa_k(\tx_k),$ and $\hat{\xi}_k = \xi_k(\hx_k)$. Also, $\bar{A}_k, \bar{G}_k, $ and $\bar{C}_k$ are the Jacobian matrices of the maps $\alpha_k, \kappa_k, $ and $\xi_k$ respectively.

Substituting \eqref{eq:rtf-lin} into \eqref{eq:rtf-2} we obtain
\begin{align}
	&\|
		x - \tx_k 
	\|^2_{\bPi_{k}}
	+
	\psi_{k}
	=
	\|
		x
	\|^2_{\bA_k^T \Omega \bA_k + \bC_k^T R_k \bC_k - \bG_k^T \bG_k} \nonumber \\
	&-
	2 x^T
	\Big(
		\left(
			\bA_k^T \Omega \bA_k
			+
			\bA_k^T \Omega
			+
			\bC_k^T R_k \bC_k 
			- 
			\bG_k^T \bG_k
		\right) \tx_k \nonumber \\
		&\;\;\;\;\;\;\;\;\;\;-
		\bA_k^T \Omega \hal_k
		+
		\bC_k^T R (y - \hxi_k)
		+
		\bG_k^T \hka_k
	\Big) \nonumber \\
	&\;+
	\|
		\tx_k
	\|^2_{\bA_k^T \Omega \bA_k + \Omega + 2 \Omega \bA_k + \bC_k^T R_k \bC_k - \bG_k^T \bG_k} \nonumber \\
	&\;-
	2 \tx_k^T
	\Big(
		\bA_k^T \Omega \hal_k
		+
		\Omega \hal_k
		-
		\bC_k^T R_k (y_k - \hxi_k)
		-
		\bG_k^T \hka_k
	\Big) \nonumber \\
\label{eq:rtf-3}
	&+
	\hal_k^T \Omega \hal_k
	+
	\| y_k \|^2_{R_k}
	-
	2 y_k^T R_k \hxi_k
	+
	\| \hxi_k \|^2_{R_k}
	-
	\hka_k^T \hka_k
	+
	\psi_{k+1}.
\end{align}
Comparing terms on the left hand side and right hand side of \eqref{eq:rtf-3}, we obtain the following recursion equations for the reverse-time filter corresponding to the time interval $[k,t]$:
\begin{align}
	&\textit{Riccati equation over the time interval}\; [k,t]: \nonumber \\
	&\bPi_{k}
	=
	\bA_k^T \Omega \bA_k 
	+ 
	\bC_k^T R_k \bC_k 
	- \bG_k^T \bG_k, \nonumber \\
\label{eq:rt-riccati}
	&\bPi_t = 0. \\
	&\textit{Filter state equation over the time interval}\; [k,t]: \nonumber \\
	&\tx_{k+1}
	=
	\tx_k
	-
	\bPi_k^{-1}
	\Lambda, \nonumber \\
\label{eq:rt-filter}
	&\tx_t = \bx_t. \\
	& \text{where} \nonumber \\
	&\Lambda 
	=
	\big(
		\bA_k^T \Omega (\tx_{k+1} - \hal_k) 
		+
		\bC_k^T R (y - \hxi_k)
		+
		\bG_k^T \bG_k	
	\big). \nonumber \\
		&\textit{Level shift scalar equation over the time interval} \; [k,t]: \nonumber \\
	&\psi_{k+1}
	=
	\psi_{k}
	-
		\|
		\tx_k
	\|^2_{\bA_k^T \Omega \bA_k + \Omega - 2 \Omega \bA_k + \bC_k^T R_k \bC_k - \bG_k^T \bG_k} \nonumber \\
	&\;\;\;
	+2 \tx_k^T
	\Big(
		\bA_k^T \Omega \hal_k
		+
		\Omega \hal_k
		-
		\bC_k^T R_k (y_k - \hxi_k)
		-
		\bG_k^T \hka_k
	\Big) \nonumber \\
	&\;\;\;
	-
	\hal_k^T \Omega \hal_k
	-
	\| y_k \|^2_{R_k}
	+
	2 y_k^T R_k \hxi_k
	-
	\| \hxi_k \|^2_{R_k}
	+
	\hka_k^T \hka_k, \nonumber \\
\label{eq:rt-scalar}
	&\psi_t = 0.
\end{align}

Hence, the robust smoothing algorithm will involve solving for two sets of Riccati differential equations, filter state equations, and level shift scalar equations. These correspond to \eqref{eq:ft-riccati}-\eqref{eq:ft-scalar} for the forward-time filter over the time interval $[0,k]$ and \eqref{eq:rt-riccati}-\eqref{eq:rt-scalar} for the reverse-time filter over the time interval $[k,t]$. 

Consider the recursion equations \eqref{eq:ft-riccati}-\eqref{eq:ft-scalar} to have solutions $\hx_k, \Pi_k, \phi_k$ over the time interval $[0,k]$, such that $\Pi_k \> 0, Q_k \> 0, R_k \> 0$. Also, if the recursion equations \eqref{eq:rt-riccati}-\eqref{eq:rt-scalar} have solutions $\tx_k, \bPi_k, \psi_k$ over the time interval $[k,t]$, such that $\bPi_k \> 0, Q_k \> 0, R_k \> 0$, then the corresponding set-valued state smoother is given by
\begin{align}
	X_{k}[\bx_0, y_{[0,t]},d]
	= 
	\Big\{
		x_{k} \in \mathbb{R}^n
		:
		\|
			&x - \hx_k
		\|^2_{\Pi_k}
		+
		\|
			x - \tx_k
		\|^2_{\bPi_k} \nonumber \\
		&\le
		d - \phi_k - \psi_k		
	\Big\}.
\label{eq:final-svse}
\end{align}
Here, the level shift scalar equations for the forward-time filter in \eqref{eq:ft-scalar} and the reverse-time filter in \eqref{eq:rt-scalar} are used to complete the definition of the set-valued state estimator \eqref{eq:final-svse} for the robust smoothing algorithm and do not contribute to the Riccati difference equations or the filter state equations for either the forward-time or the reverse-time filters.

\section{Conclusion and Future Works}
\label{sec:conclusion}
This paper derived a robust smoothing algorithm for a class of discrete-time uncertain nonlinear systems. It was assumed that the systems considered had uncertainties in modeling as well as exogenous noise sources. These uncertainties were described in terms of a sum quadratic constraint. The robust smoothing problem was solved by computing recursion equations for a forward-time and a reverse-time filter. Both filters were formulated as set-valued state estimators which were further recast into corresponding optimal control problems. These optimal control problems were represented in terms of discrete-time Hamilton-Jacobi-Bellman equations. This formulation was straightforward in the case of the reverse-time filter. However, in the case of the forward-time filter, it was necessary to employ a forward-time formulation for the corresponding Hamilton-Jacobi-Bellman equation with reverse-time system dynamics, in order to ensure the forward-time recursion equations for the filter. The final solution to the smoothing problem consisted of computing approximate solutions to the underlying dynamic programming equations. This solution led to two sets of recursion equations and involved a Riccati difference equation, a filter state equation, and a level shift scalar equation for the forward-time as well as the reverse-time filters.

As a note on future work, the robust smoothing algorithm presented in this paper can be applied to uncertain systems with unmodeled dynamics where it is needed to estimate outliers from noisy data. Such systems generally arise in the field of signal processing.

\balance
\bibliographystyle{IEEEtran}
\bibliography{jascc2013}

\begin{thebibliography}{10}
\providecommand{\url}[1]{#1}
\csname url@samestyle\endcsname
\providecommand{\newblock}{\relax}
\providecommand{\bibinfo}[2]{#2}
\providecommand{\BIBentrySTDinterwordspacing}{\spaceskip=0pt\relax}
\providecommand{\BIBentryALTinterwordstretchfactor}{4}
\providecommand{\BIBentryALTinterwordspacing}{\spaceskip=\fontdimen2\font plus
\BIBentryALTinterwordstretchfactor\fontdimen3\font minus
  \fontdimen4\font\relax}
\providecommand{\BIBforeignlanguage}[2]{{%
\expandafter\ifx\csname l@#1\endcsname\relax
\typeout{** WARNING: IEEEtran.bst: No hyphenation pattern has been}%
\typeout{** loaded for the language `#1'. Using the pattern for}%
\typeout{** the default language instead.}%
\else
\language=\csname l@#1\endcsname
\fi
#2}}
\providecommand{\BIBdecl}{\relax}
\BIBdecl

\bibitem{kalman1960}
R.~E. Kalman, ``{A New Approach to Linear Filtering and Prediction Problems},''
  \emph{Transactions of the ASME - Journal of Basic Engineering}, vol. 82(D),
  pp. 35--45, march 1960.

\bibitem{meditch1969}
S.~Meditch, \emph{{Stochastic Optimal Linear Estimation and Control}}.\hskip
  1em plus 0.5em minus 0.4em\relax New York, USA: McGraw-Hill, 1969.

\bibitem{anderson1979}
B.~D.~O. Anderson and J.~B. Moore, \emph{{Optimal Filtering}}.\hskip 1em plus
  0.5em minus 0.4em\relax Englewood Cliffs, NJ: Prentice-Hall, 1979.

\bibitem{bar-shalom1993}
Y.~Bar-Shalom and X.-R. Li, \emph{{Estimation and Tracking: Principles,
  Techniques, and Software}}.\hskip 1em plus 0.5em minus 0.4em\relax Norwood,
  MA: Artech House, 1993.

\bibitem{theodor1994}
Y.~Theodor, U.~Shaked, and C.~E. {de Souza}, ``{A Game Theory Approach to
  Robust Discrete-time {H}$^\infty$ Estimation},'' \emph{IEEE Transactions on
  Signal Processing}, vol.~42, pp. 1486--1495, June 1994.

\bibitem{savkin1995}
A.~V. Savkin and I.~R. Petersen, ``{Recursive State Estimation for Uncertain
  Systems with an Integral Quadratic Constraint},'' \emph{IEEE Transactions on
  Automatic Control}, vol.~40, no.~6, pp. 1080--1083, June 1995.

\bibitem{james1998}
M.~R. James and I.~R. Petersen, ``{Nonlinear State Estimation for Uncertain
  Systems with an Integral Constraint},'' \emph{IEEE Transactions on Acoustics,
  Speech, and Signal Processing}, vol.~46, pp. 2926--2937, 1998.

\bibitem{petersen1999}
I.~R. Petersen and A.~V. Savkin, \emph{{Robust {K}alman Filtering for Signals
  and Systems with Large Uncertainties}}, ser. {Control Engineering}.\hskip 1em
  plus 0.5em minus 0.4em\relax Birkh{\"a}user, 1999.

\bibitem{petersen2000}
I.~Petersen, V.~Ugrinovskii, and A.~Savkin, \emph{{Robust Control Design Using
  {H}$^\infty$ Methods}}, ser. {Communication and Control Engineering}.\hskip
  1em plus 0.5em minus 0.4em\relax London, UK: Springer-Verlag, 2000.

\bibitem{wang2002}
F.~Wang and V.~Balakrishnan, ``{Robust {K}alman Filters for Linear Time-Varying
  Systems with Stochastic Parametric Uncertainties},'' \emph{IEEE Transactions
  on Signal Processing}, vol.~50, no.~4, pp. 803--813, April 2002.

\bibitem{yoon2004}
M.~Yoon, V.~Ugrinovskii, and I.~Petersen, ``{Robust finite horizon minimax
  filtering for stochastic discrete time uncertain systems},'' \emph{Systems
  and Control Letters}, vol.~52, no.~2, pp. 99--112, 2004.

\bibitem{dong2006}
Z.~Dong and Z.~You, ``{Finite-horizon robust {K}alman filtering for uncertain
  discrete time-varying systems with uncertain-covariance white noises},''
  \emph{IEEE Signal Processing Letters}, vol.~13, no.~8, pp. 493--496, August
  2006.

\bibitem{souto2009}
R.~Souto, J.~Ishihara, and G.~Borges, ``{A robust extended {K}alman filter for
  discrete-time systems with uncertain dynamics, measurements and correlated
  noise},'' in \emph{{2009 American Control Conference}}, St. Louis, MO, USA,
  Jun. 2009, pp. 1888--1893.

\bibitem{ieeetac2009}
A.~G. Kallapur, I.~R. Petersen, and S.~G. Anavatti, ``{A discrete-time robust
  extended {K}alman filter for uncertain systems with sum quadratic
  constraints},'' \emph{IEEE Transactions on Automatic Control}, vol.~54,
  no.~4, pp. 850--854, Apr. 2009.

\bibitem{acc2009}
A.~Kallapur, I.~Petersen, and S.~Anavatti, ``{A discrete-time robust extended
  {K}alman filter},'' in \emph{{American Control Conference}}, St. Louis,
  Missouri, USA, Jun. 2009, pp. 3819--3823.

\bibitem{ifac2011}
A.G.Kallapur, I.G.Vladimirov, and I.R.Petersen, ``{Robust Filtering for
  Uncertain Nonlinear Systems Satisfying a Sum Quadratic Constraint},'' in
  \emph{{18th IFAC World Congress}}, Milano, Italy, 28 Aug--2 Sep 2011, pp.
  1--7, accepted.

\bibitem{acc2012}
------, ``{Robust Filtering for Continuous Time Uncertain Nonlinear Systems
  with an Integral Quadratic Constraint},'' in \emph{{American Control
  Conference}}, Montreal, Canada, Jun. 2012, pp. 4807--4812.

\bibitem{martin1983}
C.~Martin and M.~Mintz, ``{Robust filtering and prediction for linear systems
  with uncertain dynamics: A game-theoretic approach},'' \emph{IEEE
  Transactions on Automatic Control}, vol.~28, no.~9, pp. 888--896, Sep. 1983.

\bibitem{blanco2001}
E.~Blanco, P.~Neveux, and G.~Thomas, ``{H$_\infty$ smoothing for continuous
  uncertain systems},'' in \emph{{IEEE International Conference on Acoustics,
  Speech, and Signal Processing}}, vol.~6, Salt Lake City, Utah, USA, May 2001,
  pp. 3941--3944.

\bibitem{malcom2005}
W.~Malcolm, R.~Elliott, and M.~James, ``{Risk-sensitive filtering and smoothing
  for continuous-time Markov Processes},'' \emph{Information Theory, IEEE
  Transactions on}, vol.~51, no.~5, pp. 1731--1738, May 2005.

\bibitem{felsberg2006}
M.~Felsberg, P.-E. Forssen, and H.~Scharr, ``{Channel smoothing: efficient
  robust smoothing of low-level signal features},'' \emph{Pattern Analysis and
  Machine Intelligence, IEEE Transactions on}, vol.~28, no.~2, pp. 209--222,
  feb. 2006.

\bibitem{hongguo2007}
Z.~Hongguo, Z.~Huanshui, Z.~Chenghui, and S.~Xinmin, ``{Robust Filtering and
  Fixed-lag Smoothing for Linear Uncertain System with Single Delayed
  Measurement},'' in \emph{{Control Conference, 2007. CCC 2007. Chinese}}, jul
  2007, pp. 23--27.

\bibitem{garcia2010}
\BIBentryALTinterwordspacing
D.~Garcia, ``{Robust smoothing of gridded data in one and higher dimensions
  with missing values},'' \emph{Computational Statistics and Data Analysis},
  vol.~54, no.~4, pp. 1167--1178, 2010. [Online]. Available:
  \url{http://www.sciencedirect.com/science/article/pii/S0167947309003491}
\BIBentrySTDinterwordspacing

\bibitem{bertsekas1971}
D.~P. Bertsekas and I.~B. Rhodes, ``{Recursive State Estimation for a
  Set-Membership Description of Uncertainty},'' \emph{IEEE Transactions on
  Automatic Control}, vol.~16, no.~2, pp. 117-- 128, April 1971.

\bibitem{moheimani1998}
S.~O.~R. Moheimani, A.~V. Savkin, and I.~R. Petersen, ``{Robust Filtering,
  Prediction, Smoothing, and Observability of Uncertain Systems},'' \emph{IEEE
  Transactions on Circuits and Systems - I: Fundamental Theory and
  Applications}, vol.~45, no.~4, pp. 446--457, Apr. 1998.

\bibitem{acc2013}
A.~G. Kallapur and I.~R. Petersen, ``{Robust smoothing for continuous time
  uncertain nonlinear systems},'' {American Control Conference, Washington DC,
  USA, 2013, submitted}.

\bibitem{thesis2009}
A.~G. Kallapur, ``{A discrete-time robust extended Kalman filter for estimation
  of nonlinear uncertain systems},'' Ph.D. dissertation, Aerospace, Civil, and
  Mechanical Engineering, Australian Defence Force Academy, Canberra,
  Australia, 2009.

\end{thebibliography}

\end{document}